**Time-Like Constant Slope Surfaces and Space-Like Bertrand Curves in Minkowski 3-Space**

## ABSTRACT

Defining Lorentzian Sabban frame of the unit speed time-like curves on de Sitter 2-space $\mathbb{S}_1^2$ and introducing space-like height function on the unit speed time-like curves on $\mathbb{S}_1^2$, the invariants of the unit speed time-like curves on $\mathbb{S}_1^2$ and geometric properties of de Sitter evolutes of the unit speed time-like curves on $\mathbb{S}_1^2$ are studied. A relation between space-like Bertrand curves and helices is obtained. De Sitter Darboux images of space-like Bertrand curves are equal to de Sitter evolutes. The relations between time-like constant slope surfaces lying in the space-like cone and space-like Bertrand curves in Minkowski 3-space $\mathbb{R}_1^3$ are obtained.

## 1. Introduction

Huygens [1] discovered involutes (also known as evolvents) while trying to build a more accurate pendulum clock. The curve $\alpha$ is called evolute of $\tilde{\alpha}$ and the curve $\tilde{\alpha}$ is called involute of $\alpha$ if the tangent vectors are orthogonal at the corresponding points for all $s \in I \subset \mathbb{R}$. Thus the pair of $(\alpha, \tilde{\alpha})$ is called the evolute-involute pair [2].

Bertrand curves [3] are different generalization of circular helices and are particular examples of offset curves which are used in computer-aided design (CAD) and computer-aided manufacture (CAM) [4]. Izumiya and Takeuchi [5] showed that Bertrand curves can be constructed from the unit speed curves on the 2-sphere $\mathbb{S}^2$. Also they defined spherical evolutes of the unit speed curves on $\mathbb{S}^2$ and proved that these spherical evolutes are the locus of the centre of the curvatures of the unit speed curves on $\mathbb{S}^2$. In addition, they showed that the spherical Darboux images of Bertrand curves are equal to the spherical evolutes of the unit speed curves on $\mathbb{S}^2$ in Euclidean 3-space $\mathbb{R}^3$.

Babaarslan and Yayli [6] found Bertrand curves corresponding to the tangent, binormal, principal normal and Darboux images of a space curve in $\mathbb{R}^3$.

Izumiya et al. [7] defined the notion of hyperbolic evolutes and hyperbolic height function of space-like curves on the hyperbolic space $\mathbb{H}^2$. As a result, they showed that the hyperbolic evolutes are the locus of the center of geodesic curvatures in Minkowski 3-space $\mathbb{R}^3_1$.

As a generalization of the concept of helix, we can also think constant slope surfaces whose position vectors make a constant angle with the normals at each point on the surfaces. Munteanu [8] defined constant slope surfaces in $\mathbb{R}^3$ and showed that these surfaces can be constructed from unit speed curves on $\mathbb{S}^2$ such as Bertrand curves.

Babaarslan and Yayli [9] investigated the relations among Bertrand curves, spherical images and constant slope surfaces in $\mathbb{R}^3$. Also, Babaarslan et al. [10, 11] found some different characterizations of constant slope surfaces and Bertrand curves with respect to the Darboux frame. In addition, Babaarslan and Yayli [12] showed that the constant slope surfaces can be reparametized by using quaternions and homothetic motions.

Fu and Yang [13] studied space-like constant slope surfaces in $\mathbb{R}_1^3$ and classified these surfaces in the same space. Thereafter Babaarslan and Yayli [14] gave the relations among split quaternions, homothetic motions and space-like constant slope surfaces in $\mathbb{R}_1^3$.

Fu and Wang [15] gave a complete classification of time-like constant slope surfaces in $\mathbb{R}_1^3$. They showed that $S \subset \mathbb{R}_1^3$ is a time-like constant slope surface lying in the space-like cone if and only if it can be parametrized by

$$x(u,v) = u\sin\theta\left(\cos\xi(u) f(v) + \sin\xi(u) f(v) \wedge f'(v)\right), \tag{1.1}$$

where $\theta$ is a constant satisfying $\theta \in (0, \pi/2]$, $\xi(u) = \cot\theta \ln u$ and $f$ is a unit speed time-like curve on de Sitter 2-space $\mathbb{S}_1^2$. Also, there are some different characterizations of time-like constant slope surfaces in [15] since the results of them are similar to that we found [16], we only investigate the above characterization.

Constant slope surfaces have nice shapes and they are interesting in terms of differential geometry in both Euclidean and Minkowski 3-space [8].

In the present paper, we define the notions of Lorentzian Sabban frame and de Sitter evolutes of the unit time-like curves on $\mathbb{S}_1^2$ and study the invariants of the unit speed time-like curves on $\mathbb{S}_1^2$ and geometric properties of de Sitter evolutes of the unit speed time-like curves on $\mathbb{S}_1^2$. The relations among space-like Bertrand curves, helices, de Sitter Darboux images and time-like constant slope surfaces lying in the space-like cone in $\mathbb{R}_1^3$ are also given.

## 2. Basic notations, definitions and formulas

Let $\mathbb{R}_1^3$ denote the Minkowski 3-space, that is, the real vector space $\mathbb{R}^3$ endowed with the standard Lorentzian metric [17]

$$< \boldsymbol{x}, \boldsymbol{y} > = x_1 y_1 + x_2 y_2 - x_3 y_3,$$

where $\boldsymbol{x}$, $\boldsymbol{y} \in \mathbb{R}^3$. An arbitrary vector $\boldsymbol{x} \in \mathbb{R}_1^3$ is called space-like if $< \boldsymbol{x}, \boldsymbol{x} > > 0$ or $\boldsymbol{x} = 0$, time-like if $< \boldsymbol{x}, \boldsymbol{x} > < 0$ and light-like (null) if $< \boldsymbol{x}, \boldsymbol{x} > = 0$ and $\boldsymbol{x} \neq 0$. The norm (length) of a vector $\boldsymbol{x}$ is given by $\|\boldsymbol{x}\| = \sqrt{|< \boldsymbol{x}, \boldsymbol{x} >|}$.

We say that a regular curve $\alpha : I \subset \mathbb{R} \to \mathbb{R}_1^3$ is space-like, time-like or light-like if all of its velocity vector $\alpha'(t)$ is space-like, time-like or light-like, respectively. $\alpha$ is parametrized by the arc-length parameter if $\|\alpha'(s)\| = 1$ for all $s \in I \subset \mathbb{R}$. In this case, we say that $\alpha$ is a unit speed curve.

We can define the notion of the Lorentzian cross-product as follows:

$$\wedge : \mathbb{R}_1^3 \times \mathbb{R}_1^3 \to \mathbb{R}_1^3$$

$$\big((x_1, x_2, x_3), (y_1, y_2, y_3)\big) \to \begin{vmatrix} \boldsymbol{i} & \boldsymbol{j} & -\boldsymbol{k} \\ x_1 & x_2 & x_3 \\ y_1 & y_2 & y_3 \end{vmatrix} = (x_2 y_3 - x_3 y_2, x_3 y_1 - x_1 y_3, x_2 y_1 - x_1 y_2).$$

As the cross-product in Euclidean 3-space, the Lorentzian cross-product has similar algebraic and geometric properties:

**(i)** $< \boldsymbol{x} \wedge \boldsymbol{y}, \boldsymbol{z} >= \det(\boldsymbol{x}, \boldsymbol{y}, \boldsymbol{z});$

**(ii)** $\boldsymbol{x} \wedge \boldsymbol{y} = -\boldsymbol{y} \wedge \boldsymbol{x};$

**(iii)** $(\boldsymbol{x} \wedge \boldsymbol{y}) \wedge \boldsymbol{z} = - < \boldsymbol{x}, \boldsymbol{z} > \boldsymbol{y} + < \boldsymbol{y}, \boldsymbol{z} > \boldsymbol{x};$

**(iv)** $< \boldsymbol{x} \wedge \boldsymbol{y}, \boldsymbol{x} >= 0$ and $< \boldsymbol{x} \wedge \boldsymbol{y}, \boldsymbol{y} >= 0;$

**(v)** $< \boldsymbol{x} \wedge \boldsymbol{y}, \boldsymbol{x} \wedge \boldsymbol{y} >= - < \boldsymbol{x}, \boldsymbol{x} >< \boldsymbol{y}, \boldsymbol{y} > + (< \boldsymbol{x}, \boldsymbol{y} >)^2$ for all $\boldsymbol{x}, \boldsymbol{y}, \boldsymbol{z}$ in $\mathbb{R}_1^3.$

Given a unit speed curve $\alpha$ in $\mathbb{R}_1^3$, it is possible to define a Frenet frame $\{\boldsymbol{T}(s), \boldsymbol{N}(s), \boldsymbol{B}(s)\}$ associated for each point $s$. Here $\boldsymbol{T}$, $\boldsymbol{N}$ and $\boldsymbol{B}$ are the tangent, principal normal and binormal vector fields, respectively.

Let $\alpha$ be a unit speed space-like curve in $\mathbb{R}_1^3$. We assume that $\boldsymbol{T}'(s)$ is time-like. Then we have $\boldsymbol{N}(s) = \alpha''(s) / \kappa(s)$ and $\boldsymbol{B}(s) = \boldsymbol{T}(s) \wedge \boldsymbol{N}(s)$. Thus the Frenet formulae is

$$\begin{bmatrix} \boldsymbol{T}'(s) \\ \boldsymbol{N}'(s) \\ \boldsymbol{B}'(s) \end{bmatrix} = \begin{bmatrix} 0 & \kappa(s) & 0 \\ \kappa(s) & 0 & \tau(s) \\ 0 & \tau(s) & 0 \end{bmatrix} \begin{bmatrix} \boldsymbol{T}(s) \\ \boldsymbol{N}(s) \\ \boldsymbol{B}(s) \end{bmatrix},$$

where $\kappa(s) = \sqrt{- < \boldsymbol{T}'(s), \boldsymbol{T}'(s) >}$ and $\tau(s) =< \boldsymbol{N}'(s), \boldsymbol{B}(s) >$ are the curvature and torsion of the unit speed space-like curve, respectively.

The Darboux vector of this space-like curve $\alpha$ is given by $\boldsymbol{D}(s) = \tau(s)\boldsymbol{T}(s) - \kappa(s)\boldsymbol{B}(s).$ Then de Sitter Darboux image of $\alpha$ is

$$
\begin{aligned}
\boldsymbol{C} \ : \ I \ &\rightarrow \ \mathbb{S}_1^2 \\
s \ &\rightarrow \ \boldsymbol{C}(s) \ = \ \frac{\boldsymbol{D}(s)}{\left\|\boldsymbol{D}(s)\right\|}.
\end{aligned}
$$

For a general parameter $t$ of a space-like space curve $\alpha$, we can calculate the curvature and torsion as follows:

$$
\kappa(t)=\frac{\left\|\alpha'(t)\wedge\alpha''(t)\right\|}{\left\|\alpha'(t)\right\|^3}, \ \tau(t)=\frac{\det\left(\alpha'(t),\alpha''(t),\alpha'''(t)\right)}{\left\|\alpha'(t)\wedge\alpha''(t)\right\|^2}. \tag{2.1}
$$

If a time-like or space-like curve $\alpha$ in $\mathbb{R}_1^3$ is a helix, then $\tau/\kappa$ is a constant function. Conversely, let $\alpha$ be a time-like or a space-like curve with non-null principal normal vector. If $\tau/\kappa$ is constant, then $\alpha$ is a helix. Furthermore, a time-like or space-like curve $\alpha$ is a Bertrand curve if and only if there are non-zero real constants $A$, $B$ such that $A\kappa(s)+B\tau(s)=1$ for any $s\in I\subset\mathbb{R}$.

We can define de Sitter 2-space and hyperbolic space in $\mathbb{R}_1^3$, respectively as follows:

$$
\mathbb{S}_1^2=\left\{(x_1,x_2,x_3)\in\mathbb{R}_1^3 : {x_1}^2+{x_2}^2-{x_3}^2=1\right\},
$$

$$
\mathbb{H}^2=\left\{(x_1,x_2,x_3)\in\mathbb{R}_1^3 : {x_1}^2+{x_2}^2-{x_3}^2=-1\right\}.
$$

Now we define a pseudo-orthonormal frame along a time-like curve on $\mathbb{S}_1^2$. Let $f \ : \ I\rightarrow\mathbb{S}_1^2$ be a unit speed time-like curve. We denote $v$ as the arc-length parameter of $f$. Let us denote $\boldsymbol{t}(v)=f'(v)$ and we call $\boldsymbol{t}(v)$ as the unit tangent vector of $f$ at $v$. We set a vector $\boldsymbol{s}(v)=f(v)\wedge\boldsymbol{t}(v)$ and as a consequence $\boldsymbol{s}(v)\wedge\boldsymbol{t}(v)=-f(v)$, where $f$ denotes the position vector of the curve. By definition of the time-like curve $f$, we have a Lorentzian Sabban frame $\left\{f(v),\boldsymbol{t}(v),\boldsymbol{s}(v)\right\}$ along $f$. Then we have the following pseudo-spherical Frenet-Serret formulae of $f$:

$$
\begin{aligned}
&f'(v)=\boldsymbol{t}(v),\\
&\boldsymbol{t}'(v)=f(v)+\kappa_g(v)\boldsymbol{s}(v),\\
&\boldsymbol{s}'(v)=\kappa_g(v)\boldsymbol{t}(v),
\end{aligned} \tag{2.2}
$$

where $\kappa_g(v)$ is the geodesic curvature of the unit speed time-like curve $f$ on $\mathbb{S}_1^2$ which is given by $\kappa_g(v)=\det\left(f(v),\boldsymbol{t}(v),\boldsymbol{t}'(v)\right)$.

Also we define a curve on $\mathbb{S}_1^2$ as follows:

$$d_f(v) = \frac{\kappa_g(v) f(v) - \boldsymbol{s}(v)}{\sqrt{\kappa_g^2(v)+1}}.$$

We call $d_f$ as the de Sitter evolute of the time-like curve $f$.

## 3. Space-like height function of unit speed time-like curves on $\mathbb{S}_1^2$

We introduce a function on a time-like curve $f : I \to \mathbb{S}_1^2$ [7]. Now we define a function $H^S : I \times \mathbb{S}_1^2 \to \mathbb{R}$ by $H^S(v,u) = < f(v), u >$. We call $H^S$ as the space-like height function of the time-like curve $f$ and denoted by $(h_u^S)(v) = H^S(v,u).$

Thus, we have the following proposition:

**Proposition 3.1.** Let $f : I \to \mathbb{S}_1^2$ be a unit speed time-like curve. For any $(v,u) \in I \times \mathbb{S}_1^2$ :

**(a)** $(h_u^S)'(v) = 0$ if and only if $u \in \mathrm{span}\{f(v), \boldsymbol{s}(v)\},$

**(b)** $(h_u^S)'(v) = (h_u^S)''(v) = 0$ if and only if $u = \pm\left(\kappa_g(v) f(v) - \boldsymbol{s}(v)\right) \big/ \sqrt{\kappa_g^2(v)+1}.$

**Proof.** By using Eq. (2.2), we have

**(i)** $(h_u^S)'(v) = < \boldsymbol{t}(v), u >,$

**(ii)** $(h_u^S)''(v) = < f(v) + \kappa_g(v)\boldsymbol{s}(v), u > .$

The assertion (a) can be found from the formula (i). By using assertion (a), there exist $\lambda, \mu \in \mathbb{R}$ such that $u = \lambda f(v) + \mu \boldsymbol{s}(v).$ From formula (ii), we obtain

$$\begin{aligned} 0 &= < f(v) + \kappa_g(v)\boldsymbol{s}(v), \lambda f(v) + \mu \boldsymbol{s}(v) > \\ &= \lambda < f(v), f(v) > + \mu\kappa_g(v) < \boldsymbol{s}(v), \boldsymbol{s}(v) > \\ &= \lambda + \mu\kappa_g(v). \end{aligned}$$

Thus we get $u = -\mu\left(\kappa_g(v) f(v) - \boldsymbol{s}(v)\right).$ Since $< u,u > = 1,$ we have

$$\mu = \mp \frac{1}{\sqrt{\kappa_g^2(v)+1}}.$$

As a result, we have

$$u = \pm \frac{1}{\sqrt{\kappa_g^2(v)+1}} \left( \kappa_g(v) f(v) - \boldsymbol{s}(v) \right).$$

Conversely, substituting $u$ into (i) and (ii), respectively, we find $(h_u^S)'(v) = (h_u^S)''(v) = 0.$ This completes the proof.

## 4. Spherical invariants of unit speed time-like curves on $\mathbb{S}_1^2$

We study the geometric properties of de Sitter evolutes of the unit speed time-like curves on $\mathbb{S}_1^2$ [7]. For any $r \in \mathbb{R}$ and $u_0 \in \mathbb{S}_1^2$, we denote $PS^1(u_0, r) = \{u \in \mathbb{S}_1^2 : < u, u_0 >= r\}$. We call $PS^1(u_0, r)$ as a pseudo-circle whose center is $u_0$ on $\mathbb{S}_1^2$.

Then we have the following proposition:

**Proposition 4.1.** Let $f : I \to \mathbb{S}_1^2$ be a unit speed time-like curve. Then $\kappa_g'(v) \equiv 0$ if and only if $u_0 = \pm \left( \kappa_g(v) f(v) - s(v) \right) \Big/ \sqrt{\kappa_g^2(v)+1}$ are constant vectors. Under this condition, $f$ is a part of a pseudo-circle whose center is $u_0$ on $\mathbb{S}_1^2$.

**Proof.** Let we denote

$$P_{\pm}(v) = \pm u_0 = \pm \frac{1}{\sqrt{\kappa_g^2(v)+1}} \left( \kappa_g(v) f(v) - s(v) \right).$$

Taking the derivative of this equation with respect to $v$, we have

$$P_{\pm}'(v) = \pm \kappa_g'(v) \frac{\left( f(v) + \kappa_g(v) s(v) \right)}{\left( \kappa_g^2(v)+1 \right)^{3/2}}.$$

Thus $P_{\pm}'(v) \equiv 0$ if and only if $\kappa_g'(v) \equiv 0.$

Under this condition, we put $r = \pm \kappa_g(v) \Big/ \sqrt{\kappa_g^2(v)+1}$ and $u_0 = \pm \left( \kappa_g(v) f(v) - s(v) \right) \Big/ \sqrt{\kappa_g^2(v)+1}.$ Then $f(v)$ is a part of the pseudo-circle $PS^1(u_0, r).$ This completes the proof.

Let $f : I \to \mathbb{S}_1^2$ be a unit time-like curve. For any $v_0 \in I,$ we consider the pseudo-circle $PS^1(u_0, r_0^{\pm}),$ where $u_0 = d_f(v_0)$ and $r_0 = \kappa_g(v_0) \Big/ \sqrt{\kappa_g^2(v_0)+1}.$

Thus we have the following proposition:

**Proposition 4.2.** Under the above notations, $f$ and $PS^1(u_0, r_0)$ have at least a 3-point contact at $f(v_0)$.

**Proof.** Proposition 3.1 (b) says that $f$ and $PS^1(u_0, r_0)$ have at least a 3-point contact at $f(v_0)$. Thus, the proof is completed.

**Remark 4.3.** We call $PS^1(u_0, r_0)$ in proposition 4.2 as the pseudo-circle of geodesic curvature and its center $u_0$ is called as the center of geodesic curvature. As a result, the de Sitter evolute is the locus of the center of geodesic curvature.

## 5. Time-like constant slope surfaces lying in the space-like cone and space-like Bertrand curves

In this section, we give the relations among space-like Bertrand curves, helices, de Sitter Darboux images and time-like constant slope surfaces lying in the space-like cone in $\mathbb{R}^3_1$.

Now we can express the following lemma:

**Lemma 5.1.** Let $f : I \to \mathbb{S}^2_1$ be a unit speed time-like curve. Then

$$\tilde{\gamma}(v) = a\int_0^v f(t)dt + a\tan\xi\int_0^v f(t)\wedge f'(t)dt \tag{5.1}$$

is a space-like Bertrand curve, where $a$ and $\xi = \xi(u) = \cot\theta\ln u$ are constant numbers, and $\theta$ is a constant satisfying $\theta\in(0, \pi/2]$. Moreover, all space-like Bertrand curves can be constructed by using this method.

**Proof.** We now calculate the curvature and the torsion of $\tilde{\gamma}(v)$. Taking the derivatives of Eq. (5.1) three times with respect to $v$, we have

$$\tilde{\gamma}'(v) = a\left(f(v) + \tan\xi\, \boldsymbol{s}(v)\right),$$
$$\tilde{\gamma}''(v) = a\left(1 + \tan\xi\, \kappa_g(v)\right)\boldsymbol{t}(v),$$
$$\tilde{\gamma}'''(v) = a\left(\left(1 + \tan\xi\, \kappa_g(v)\right) f(v) + \tan\xi\, \kappa_g'(v)\boldsymbol{t}(v) + \left(\kappa_g(v) + \tan\xi\, \kappa_g^2(v)\right)\boldsymbol{s}(v)\right).$$

Therefore, by using Eq. (2.1), we obtain $\kappa(v)$ and $\tau(v)$ as follows:

$$\kappa(v) = \varepsilon \frac{\cos^2\xi\left(1 + \tan\xi\, \kappa_g(v)\right)}{a} \text{ and } \tau(v) = \frac{\cos^2\xi\left(\kappa_g(v) - \tan\xi\right)}{a}, \tag{5.2}$$

where $\varepsilon = \pm 1$. It follows from these formulae that $a\left(\varepsilon\kappa(v) - \tan\xi\, \tau(v)\right) = 1$, thus $\tilde{\gamma}(v)$ is a Bertrand curve. Also, since $< \tilde{\gamma}'(v), \tilde{\gamma}'(v) >= a^2 / \cos^2\xi > 0$, $\tilde{\gamma}(v)$ is a space-like Bertrand curve.

Conversely, let $\tilde{\gamma}(s)$ be a space-like Bertrand curve. Thus, there exist real constants $A$, $B$ different from zero such that $A\kappa(s) + B\tau(s) = 1$. Here we put $A = a$ and $B = -a\tan\xi$. Assume that $a > 0$ and choose $\varepsilon = \pm 1$ with $\varepsilon \cos\xi / a > 0$.

Let us consider the Frenet frame $\{\boldsymbol{T}(s), \boldsymbol{N}(s), \boldsymbol{B}(s)\}$ for the space-like Bertrand curve $\tilde{\gamma}(s)$. In this trihedron $\boldsymbol{T}(s)$ and $\boldsymbol{B}(s)$ are space-like vectors, $\boldsymbol{N}(s)$ is a time-like vector. For these vectors, we have

$$\boldsymbol{T}(s) \wedge \boldsymbol{N}(s) = \boldsymbol{B}(s) \text{ and } \boldsymbol{B}(s) \wedge \boldsymbol{N}(s) = -\boldsymbol{T}(s).$$

Now we define a time-like curve on $\mathbb{S}_1^2$ as

$$f(s) = \varepsilon\left(\cos\xi\, \boldsymbol{T}(s) - \sin\xi\, \boldsymbol{B}(s)\right).$$

Thus we have

$$f'(s) = \varepsilon\cos\xi\left(\kappa(s) - \tan\xi\, \tau(s)\right)\boldsymbol{N}(s) = \frac{\varepsilon}{a}\cos\xi\, \boldsymbol{N}(s).$$

Let $v$ be the arc-length parameter of $f$, then we have $dv / ds = \varepsilon\cos\xi / a$. Moreover, we get

$$af(s)\frac{dv}{ds} = \cos\xi\left(\cos\xi\, \boldsymbol{T}(s) - \sin\xi\, \boldsymbol{B}(s)\right) \tag{5.3}$$

and

$$a\tan\xi\, f(s) \wedge \frac{df}{dv}\frac{dv}{ds} = \sin\xi\left(\cos\xi\, \boldsymbol{B}(s) + \sin\xi\, \boldsymbol{T}(s)\right). \tag{5.4}$$

By using Eq. (5.3) and (5.4), we obtain

$$\begin{aligned} a\int_0^v f(t)dt + a\tan\xi\int_0^v f(t)\wedge f'(t)dt &= \int_{s_0}^{s}\cos\xi\left(\cos\xi\boldsymbol{T}(t)-\sin\xi\boldsymbol{B}(t)\right)dt \\ &\quad +\int_{s_0}^{s}\sin\xi\left(\cos\xi\boldsymbol{B}(t)+\sin\xi\boldsymbol{T}(t)\right)dt \\ &= \int_{s_0}^{s}\boldsymbol{T}(t)dt = \tilde{\gamma}(s). \end{aligned}$$

This completes the proof.

As a consequence of this lemma, we can give a relation between space-like Bertrand curves and helices.

**Corollary 5.2.** The unit speed time-like curve $f$ on $\mathbb{S}_1^2$ is a part of a pseudo-circle if and only if the corresponding space-like Bertrand curve is a helix.

**Proof.** Taking the derivative of Eq. (5.2) with respect to $v,$ we have

$$\kappa'(v)=\varepsilon\frac{\sin 2\xi\kappa_g'(v)}{2a} \text{ and } \tau'(v)=\frac{\cos^2\xi\kappa_g'(v)}{a}.$$

From proposition 4.1, the unit speed time-like curve $f$ on $\mathbb{S}_1^2$ is a part of a pseudo-circle if and only if $\kappa_g'(v)\equiv 0.$ This condition is equivalent to the condition that both $\kappa(v)$ and $\tau(v)$ are constants. The proof is completed.

Then we have the following proposition:

**Proposition 5.3.** Let $f : I\to\mathbb{S}_1^2$ be a unit speed time-like curve and $\tilde{\gamma}: I\to\mathbb{R}_1^3$ be a space-like Bertrand curve corresponding to $f.$ Then the de Sitter Darboux image of $\tilde{\gamma}$ is equal to the de Sitter evolute of $f.$

**Proof.** From Eq. (5.2), we have

$$\kappa(v)=\varepsilon\frac{\cos^2\xi\left(1+\tan\xi\kappa_g(v)\right)}{a} \text{ and } \tau(v)=\frac{\cos^2\xi\left(\kappa_g(v)-\tan\xi\right)}{a}.$$

For the space-like curve $\tilde{\gamma},$ we obtain

$$\boldsymbol{T}(v)=a\left(f(v)+\tan\xi\boldsymbol{s}(v)\right)\frac{dv}{ds} \text{ and } \boldsymbol{N}(v)=\varepsilon\boldsymbol{t}(v).$$

Then we get

$$\boldsymbol{B}(v)=\boldsymbol{T}(v)\wedge \boldsymbol{N}(v)=\varepsilon a\frac{dv}{ds}\big(\boldsymbol{s}(v)-\tan\xi f(v)\big).$$

We can easily show that

$$\boldsymbol{D}(v)=\tau(v)\boldsymbol{T}(v)-\kappa(v)\boldsymbol{B}(v)=\frac{dv}{ds}\big(\kappa_g(v)f(v)-\boldsymbol{s}(v)\big).$$

Therefore we have $\boldsymbol{C}(v)=\boldsymbol{D}(v)/\|\boldsymbol{D}(v)\|=d_f(v).$ This completes the proof.

We have the following theorem:

**Theorem 5.4.** Let $f: I\to\mathbb{S}_1^2$ be a unit speed time-like curve and $\tilde{\gamma}: I\to\mathbb{R}_1^3$ be a space-like Bertrand curve corresponding to $f.$ Then $\tilde{\gamma}'(v)$ lies on the time-like constant slope surface $x(u,v)$ lying in the space-like cone.

**Proof.** Taking the derivative of Eq. (5.1) with respect to $v,$ we obtain

$$\tilde{\gamma}'(v)=af(v)+a\tan\xi f(v)\wedge f'(v).$$

In this equation, we can take $a$ as $a=u\sin\theta\cos\xi$ and so $a\tan\xi=u\sin\theta\sin\xi,$ where $u,$ $\theta$ are constants. Thus by Eq. (1.1), $\tilde{\gamma}'(v)$ is $v-$parameter curve of time-like constant slope surface $x(u,v)$ lying in the space-like cone and $\tilde{\gamma}'(v)$ lies on it. This completes the proof.

We now state the relation between time-like constant slope space-like surfaces lying in the space-like cone and space-like Bertrand curves.

**Theorem 5.5.** Let $x: S\to\mathbb{R}_1^3$ be a time-like constant slope surface immersed in $\mathbb{R}_1^3$ and $x$ lies in the space-like cone. If $x(v)$ is $v-$parameter curve of time-like constant slope surface $x(u,v)$ lying in the space-like cone, then $\int_0^v x(v)dv$ is a space-like Bertrand curve.

**Proof.** From Eq. (1.1), we have

$$x(v)=u\sin\theta\cos\xi f(v)+u\sin\theta\sin\xi f(v)\wedge f'(v)$$

for $u=$constant, where $\xi=\xi(u)=\cot\theta\ln u.$ By integrating $x(v),$ we have the equation as

$$\int_0^v x(v)dv=u\sin\theta\cos\xi\int_0^v f(v)dv+u\sin\theta\sin\xi\int_0^v f(v)\wedge f'(v)dv.$$

Since the coefficients of $f(v)$ and $f(v)\wedge f'(v)$ are constants, here we can take $u\sin\theta\cos\xi$ as $u\sin\theta\cos\xi = a$ and so $u\sin\theta\sin\xi = a\tan\xi$. Thus we obtain

$$\int_0^v x(v)dv = a\int_0^v f(v)dv + a\tan\xi\int_0^v f(v)\wedge f'(v)dv.$$

From lemma 5.1, $\int_0^v x(v)dv$ is a space-like Bertrand curve. This completes the proof.

Now we give an example of time-like constant slope surfaces and space-like Bertrand curves and draw the corresponding pictures via Mathematica.

**Example 5.6.** By using Eq. (1.1), we may choose the unit speed time-like curve as $f(v) = (\cosh v, 0, \sinh v)$ on $\mathbb{S}_1^2$. Then we have $f(v)\wedge f'(v) = (0,-1,0)$. Thus, the time-like constant slope surface lying in the space-like cone is

$$x(u,v) = u\sin\theta\left(\cos(\cot\theta\ln u)\cosh v, -\sin(\cot\theta\ln u), \cos(\cot\theta\ln u)\sinh v\right).$$

For $\theta = \pi/4$, the picture of this surface is given in Fig. 1. Thus, for $u = e$, the space-like Bertrand curve is

$$\int_0^v x(v)dv = \frac{\sqrt{2}}{2}e\left(\cos(1)\sinh v, -\sin(1)v, \cos(1)(\cosh v - 1)\right).$$

Since the time-like curve $f(v)$ is a part of a pseudo-circle, from corollary 5.2, this space-like Bertrand curve is a helix. The picture of this curve is given in Fig. 2.

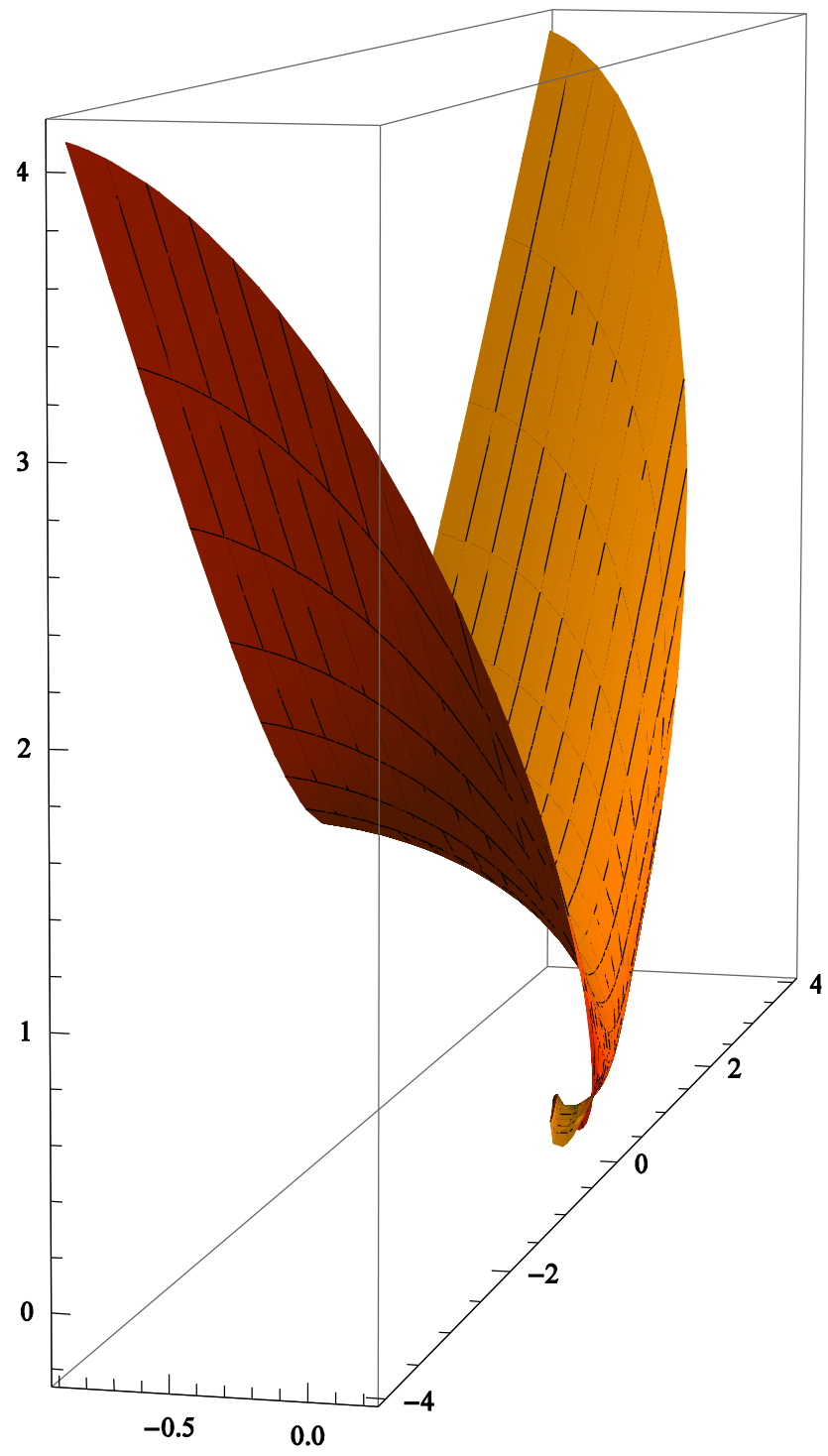

**Fig. 1.** Time-like constant slope surface lying in the space-like cone, $f(v) = (\cosh v, 0, \sinh v)$, $\theta = \pi / 4$

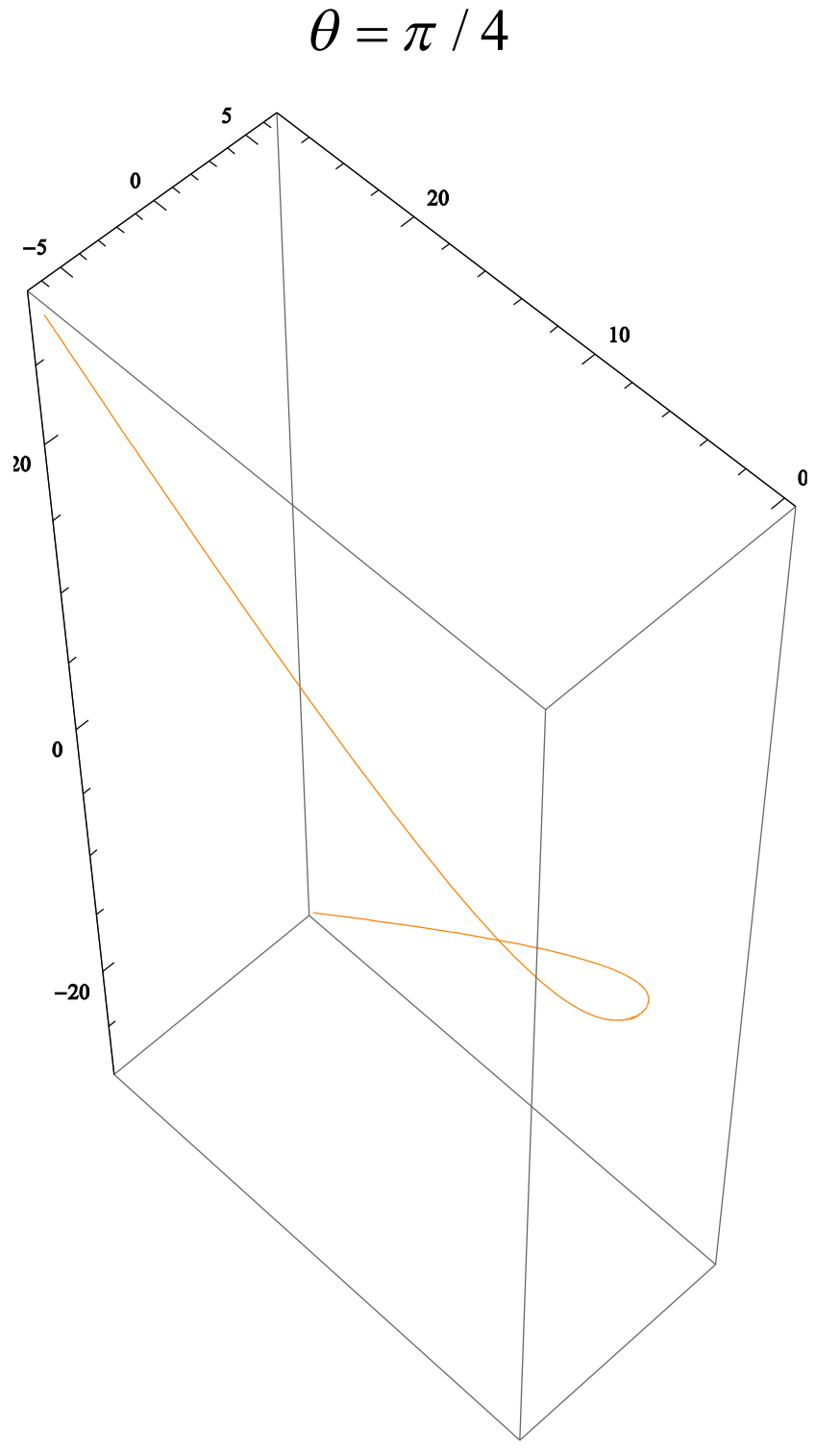

**Fig. 2.** Space-like Bertrand curve, $\theta = \pi / 4$, $u = e$